\documentclass[a4paper,twosided,11pt]{article}

\usepackage[utf8]{inputenc} 
\usepackage[T1]{fontenc} 
\usepackage[a4paper]{geometry} 
\usepackage{amsmath,amssymb,mathtools,amsthm}
\usepackage{mathrsfs}  
\usepackage{graphicx} 
\usepackage{xcolor}   
\usepackage[english]{babel} 
\usepackage{bbm}
\usepackage[title]{appendix}
\usepackage{caption} 
\usepackage{subcaption} 
\usepackage{comment}
\usepackage{cite}
\usepackage{empheq}
\usepackage{float}
\usepackage{booktabs}
\usepackage{varwidth}
\usepackage{enumitem} 
\usepackage[colorlinks=true]{hyperref} 

\newcommand \commentout[1] {}

\newcommand{\R}{\mathbb{R}}
\newcommand{\N}{\mathbb{N}}

\newcommand{\T}{\mathbb{T}}
\newcommand{\C}{\mathbb{C}}

\newcommand{\al}{\alpha}
\newcommand{\be}{\beta}

\newcommand{\eps}{\varepsilon}

\newcommand{\la}{\lambda}

\newcommand{\Chi}{{\bf \raise 2pt \hbox{$\chi$}}}

\newcommand{\p}{\partial}

\newcommand*{\dd}{\mathop{\kern0pt\mathrm{d}}{}}
\newcommand*{\DD}{\mathop{\kern0pt\mathrm{D}}{}}

\newcommand{\av}[1]{\left\langle #1 \right\rangle}
\newcommand{\avmu}[1]{\left\langle #1 \right\rangle_{\mu}}

\DeclareMathOperator*{\Cov}{\operatorname{Cov}}

\theoremstyle{plain}
\newtheorem*{thm*}{Theorem}
\newtheorem{thm}{Theorem}[section]
\newtheorem{lemma}[thm]{Lemma}

\newtheorem{proposition}[thm]{Proposition}
\newtheorem{corollary}[thm]{Corollary}

\theoremstyle{remark}
\newtheorem{remark}[thm]{\bf Remark}
\newtheorem{definition}[thm]{\bf Definition}

\newcommand{\beq}{\begin{equation}}
\newcommand{\eeq}{\end{equation}}
\newcommand{\bea}{\begin{array}{rl}}
\newcommand{\eea}{\end{array}}
\newcommand{\bepa}{\left\{ \begin{array}{l}}
\newcommand{\eepa}{\end{array}\right.}
\newcommand{\diff}{\mathop{} \mathrm{d}}

\newcommand{\weakto}{\rightharpoonup}

\numberwithin{equation}{section}

\title{Global existence for a system without self-diffusion and different mobilities}
\author{
    Charles Elbar%
     \thanks{Universite Claude Bernard Lyon 1, ICJ UMR5208, CNRS, Ecole Centrale de Lyon, INSA Lyon, Université Jean Monnet, 69622
Villeurbanne, France. Email: elbar@math.univ-lyon1.fr} %
}
\date{}

\begin{document}
\maketitle
    
\begin{abstract}
We study a one-dimensional cross-diffusion system for two interacting populations on the torus, with a linear pressure law and different mobilities.  For arbitrary bounded non-negative initial data, we show that any good approximation scheme, yields existence of global weak solutions. More precisely, we introduce a notion of \textit{admissible approximation sequence} and show that any such sequence admits a subsequence converging to a weak solution of the system. The strategy relies on entropy estimates and the div--curl lemma, in the framework of Young measures.
\end{abstract}
\vskip .7cm

\noindent{\makebox[1in]\hrulefill}\newline

\section{Introduction}

We consider the cross-diffusion system
\begin{equation}
\label{eq:MNsystem}
\begin{cases}
\partial_t m-\partial_x \left(m \partial_x(m+n)\right)=0,\\
\partial_t n-\nu \partial_x \left(n \partial_x(m+n)\right)=0,
\end{cases}
\quad (t,x)\in(0,\infty)\times\mathbb T,
\end{equation}
posed on the flat torus $\mathbb T=\mathbb R/\mathbb Z$, with parameter $\nu>0$, $\nu\neq1$, and nonnegative initial data $m_0,\, n_0$. The unknowns $m$ and $n$ are nonnegative densities, and we denote by
$$
\rho:=m+n
$$
their total density. The two species are transported by the same pressure gradient $\partial_x\rho$, but with different mobilities: the first species moves with mobility $1$, whereas the second moves with mobility $\nu$. 

Systems of this type appear in several contexts where the motion of each population is driven by crowding effects, rather than by independent Fick's type of diffusion. In ecology and population dynamics, this point of view goes back to the work of Gurtin and Pipkin and the work of Busenberg and Travis ~\cite{GurtinPipkin1984,BusenbergTravis1983}. They also appear in models for cell motion, tissue growth, crowd motion~\cite{MR4712820, MR4880211,GwiazdaPerthameSwierczewska2019,lorenzi_interfaces_2016}. Similar systems have also been derived as local limits of nonlocal or particle models, \cite{ChenDausJungel2019,ChenDausHolzingerJungel2021,DoumicHechtPerthamePeurichard2024}.

\eqref{eq:MNsystem} can be rewritten as
$$
\partial_t \begin{pmatrix}
m\\ n
\end{pmatrix}
- \partial_x \left(
\begin{pmatrix} m&m\\
\nu n&\nu n
\end{pmatrix} \partial_x
\begin{pmatrix}
m\\ n
\end{pmatrix}
\right)=0.
$$
The diffusion matrix has rank one. In particular, only one diffusive mode is present, whereas the second mode is transported but not regularized. The classical entropy
$$
\mathcal H(m,n):=
\int_{\mathbb T}\left(h(m)+\frac1\nu h(n)\right) \mathrm dx,
\quad
h(z)=z\log z-z,
$$
formally satisfies
$$
\frac{\mathrm d}{\mathrm dt}\mathcal H(m,n)
+\int_{\mathbb T}|\partial_x(m+n)|^2 \mathrm dx=0.
$$
Thus, the entropy controls $\partial_x\rho$, but not $\partial_xm$ and $\partial_xn$ separately. This is the main source of difficulty for compactness arguments and to prove global existence of the system.

When $\nu=1$, the total density solves
$$
\partial_t\rho-\partial_x(\rho \partial_x\rho)=0.
$$
Then, one can obtain strong compactness of a sequence $\partial_x\rho^k$, which is enough to pass to the limit in the nonlinear fluxes
$$
m^k \partial_x\rho^k,
\quad
n^k \partial_x\rho^k.
$$
The case $\nu\neq1$ is substantially harder. Only weak compactness of $\partial_x\rho^k$ is available anymore. Even if one has strong compactness of $\rho^k$, this does not help to pass to the limit in the nonlinear fluxes. The difficulty is therefore to understand the oscillations between the densities of the species and the pressure gradient.

\subsection{Previous results on this model and on Busenberg--Travis systems}

In \cite{MR821681}, Bertsch, Gurtin, Hilhorst and Peletier studied System~\eqref{eq:MNsystem} and proved the existence of a solution in which the two species are segregated for all time, provided it is initially the case. This result is independent of $\nu$. This is one of the main results on the topic, and it already shows that the rank one matrix can lead to sharp interfaces and non-smooth solutions. In another work \cite{MR887658}, Bertsch, Gurtin and Hilhorst treated the special case $\nu=1$ with equal mobility, and proved existence of solutions for general initial data. In~\cite{lorenzi_interfaces_2016}, Lorenzi, Lorz and Perthame showed the existence of traveling waves for System~\eqref{eq:MNsystem} with a source term. Druet, Hopf and Jüngel proved local well-posedness for $H^s$ initial data \cite{DruetHopfJungel2023}. Hopf and Jüngel~\cite{HopfJungel2025} then developed a measure-valued framework for this hyperbolic--parabolic class of systems, proved weak--strong uniqueness, and showed convergence of a finite-volume approximation scheme to measure-valued solutions. These results somehow already suggest that, measure-valued notions are important to understand the defect of compactness.

Other works considered variants or regularizations of the Busenberg--Travis system. Jüngel, Vetter and Zurek studied a nonlocal regularization of a generalized Busenberg--Travis system \cite{2024arXiv240701123J}. In another direction, Carrillo, Chen, Du and Jüngel obtained a fluid-relaxation approximation of the Busenberg--Travis system on the torus, using compactness and relative-entropy arguments \cite{CarrilloChenDuJungel2025}.

\subsection{A useful change of variables}

To see the source of oscillation, it is natural to introduce, besides $\rho$, the scalar variable
$$
a=\frac{m+\nu n}{m+n}=\frac{m+\nu n}{\rho},
$$
which we call the \emph{activity variable}. Thus the existence for $(m,n)$ can be reformulated for $(\rho,a)$. This is useful because the total density satisfies
$$
\partial_t\rho-\partial_x \left(a\rho \partial_x\rho\right)=0,
$$
so that $a$ is a mobility for the only diffusive variable $\rho$. With standard approximations, one expects weak compactness of $\partial_x \rho^k$, strong compactness of $\rho^k$, but only weak compactness of $a^k$. If one is able to show that the weak limit of $a^k \rho^k\partial_x \rho^k$ is $a \rho \partial_x \rho$ then the previous equation holds in the limit.  Then one can pass to the limit in System~\eqref{eq:MNsystem} as for instance the flux $m^k\partial_x\rho^k$ can be rewritten as 
$$
m^k\partial_x\rho^k = \frac{\nu}{\nu-1}\rho^k\partial_x\rho^k - \frac{1}{\nu-1}a^k\rho^k\partial_x\rho^k.
$$

As $\rho^k$ converges strongly we see that the only difficulty is to prove that $a^k\partial_x\rho^k$ converges weakly to $a\partial_x\rho$ which is the main result of this paper.  

The defect, which is due to oscillations and not concentrations as $a^k$ and $\partial_x \rho^k$ remain bounded (respectively in $L^\infty$ and $L^2$), is therefore concentrated in the vector
$$
(a^k,\xi^k), \quad \xi^k=\partial_x\rho^k.
$$
The question is then whether the oscillations of $a^k$ can be coupled strongly enough to those of $\xi^k$ to identify the weak limits.

\subsection{Young measures}

Young measures are a natural tool to study these oscillations. This framework is necessary since in the usual approximation schemes we  cannot pass to the limit in the terms of the form $f(z^k)$ where $f$ is nonlinear and $\{z^k\}$ is only a weak star convergent sequence. The idea of Young measures is to embed the problem in a larger space and gain linearity. We write $f(z^k(y))=\langle f,\delta_{z^k(y)}\rangle$, and if $f\in C_{0}(\R^n)$, using the duality $(L^{1}(Q;C_{0}(\R^n)))^{*}=L^{\infty}_{w}(Q;\mathcal{M}(\R^n))$, Banach-Alaoglu theorem and weak-star continuity of linear operators, we can pass to the limit. The same is true if $\{f(z_j)\}$ is weakly compact in $L^1(Q)$ and $\{z_j\}$ does not grow too fast as the following theorem states:
\begin{thm}[Fundamental Theorem of Young Measures]\label{thm:fund_theorem_YM}
Let $Q\subset\R^{d}$ be a measurable set and let $z^k:Q\to \R^{n}$ be measurable functions such that 
\begin{equation*}
\sup_{k\in\mathbb{N}}\int_{Q}g(|z^k(y)|)\mathrm{d}y<+\infty   
\end{equation*}
for some continuous, nondecreasing function $g:[0,+\infty)\to [0,+\infty)$ with $\lim_{t\to+\infty}g(t)=+\infty$. Then, there exists a subsequence (not relabeled) and a {weak star} measurable family of probability measures $\nu= \{\nu_{y}\}_{y\in Q}$ with the property that whenever the sequence $\{\psi(y, z^k (y))\}_{k\in\mathbb{N}}$ is weakly compact in $L^{1}(Q)$ for a Carathéodory function (measurable in the first and continuous in the second argument) $\psi: Q \times \R^{n}\to \R$, we have
\begin{equation*}
\psi(y,z^k(y))\rightharpoonup\int_{\R^{n}}\psi(y,\lambda) \mathrm{d}\nu_{y}(\lambda)\quad \text{in $L^{1}(Q)$}.     
\end{equation*}
We say that the sequence $\{z^k\}_{j\in\mathbb{N}}$ generates the Young measure $\{\nu_{y}\}_{y\in Q}$.
\end{thm}

Note that in the above theorem, we require that the sequence $\{\psi(y,z^k(y))\}$ is weakly compact in $L^{1}(Q)$. This prevents the concentration effect to appear (think about the family of standard mollifiers). Also, we note the important fact that if one is able to prove that the Young measures are Dirac masses, then the limit of the nonlinearities can be identified, proving therefore strong convergence of the sequence.

In our case, we encode the possible oscillations of the pair $(a^k,\xi^k)$ by a parametrized Young measure
$$
\lambda_{t,x}\in\mathcal P(I\times\mathbb R).
$$

The use of Young measures in nonlinear PDE goes back to Tartar's compensated compactness and to Ball's abstract Young-measure theorem, and it has become standard whenever nonlinear quantities must be identified without strong convergence~\cite{Tartar1979,Ball1989}. In systems without self-diffusion, the measure-valued viewpoint has recently become especially relevant. On the one hand, Hopf and Jüngel used dissipative measure-valued solutions to describe global limits in~\cite{HopfJungel2025}. On the other hand, several recent papers show that the oscillations often have a particular structure.

In particular, Mészáros and Parker, and later the author with Santambrogio, studied one-dimensional two-species systems with diffusion acting only on the total density and with different external drifts. Although those works did not use Young measures, they proved global existence of solutions~\cite{2025arXiv251007937E,MeszarosParker2025,2026arXiv260318770M}, observing that a particular ratio variable does not oscillate and remains of bounded variation. Very recently, Skrzeczkowski considered one-dimensional cross-diffusion systems with independent drifts and any pressure with power law on the total density. He showed existence using Young measures~\cite{Skrzeczkowski2026}. These works are not about \eqref{eq:MNsystem} itself, but they support the idea that for one-dimensional systems with rank one diffusion, the right target is not strong compactness of each species, but rather compactness of an appropriately chosen scalar quantity.

\subsection{Notations and main results}

We work on the torus $\mathbb T=\mathbb R/\mathbb Z$. The total density and activity variables are
$$
\rho:=m+n, \quad a:=\frac{m+\nu n}{m+n},
$$
so that

\begin{equation}
\label{eq:mnfromrhoa}
 m=\frac{\nu-a}{\nu-1} \rho,  \quad n=\frac{a-1}{\nu-1} \rho.
\end{equation}

We write
$$
\xi:=\partial_x\rho.
$$
The activity takes values in
$$
I=[\alpha,\beta],
\quad \alpha:=\min\{1,\nu\}, \quad \beta:=\max\{1,\nu\}.
$$

For approximation sequences we write
$$
\rho^k=m^k+n^k,
\quad a^k=\frac{m^k+\nu n^k}{\rho^k},
\quad \xi^k=\partial_x\rho^k.
$$
The Young measure generated by $(a^k,\xi^k)$ is denoted by
$$
\lambda_{t,x}\in\mathcal P(I\times\mathbb R).
$$

We note for continuous bounded function $f:I\times\R\to\R$
$$
f(a^k,\xi^k)
\rightharpoonup
\int_{I\times\R} f(\lambda_1,\lambda_2)\,\diff\lambda_{t,x}(\lambda_1,\lambda_2)
\quad \text{weakly in } L^1\bigl((0,T)\times \T\bigr).
$$
Later, after disintegration we use the notations

$$
\av{f(a,\xi)}:=\int_{I\times\R} f(\lambda_1,\lambda_2)\,\diff\lambda_{t,x}(\lambda_1,\lambda_2),
\quad
\avmu{g(a)}:=\int_I g(\lambda_1)\,\diff\mu_{t,x}(\lambda_1).
$$

Rather than building a specific approximation scheme, we just focus on the compactness needed. This leads to the notion of an admissible approximation sequence in Definition~\ref{def:admissible}. The theorem below says that any sequence satisfying these properties converges, up to a subsequence, to a weak solution.

\begin{definition}
\label{def:weak-solution}
Let $T>0$. A pair of nonnegative functions $(m,n)$ on $(0,T)\times\T$ is a \emph{weak solution} of \eqref{eq:MNsystem} with initial data $(m_0,n_0)$ if the following hold:
\begin{enumerate}[label=(\roman*)]
\item $m,n\in L^\infty((0,T)\times\T)$, and $\rho=m+n\in L^2(0,T;H^1(\T))$;
\item for every test function $\varphi\in C_c^\infty([0,T)\times\T)$,
\begin{equation}
 -\int_0^T  \int_{\T} m\,\p_t\varphi\,\diff x\diff t
 -\int_{\T}m_0(x)\varphi(0,x)\,\diff x
 +\int_0^T  \int_{\T} m\partial_x\rho\,\p_x\varphi\,\diff x\diff t=0,
\end{equation}
and
\begin{equation}
 -\int_0^T  \int_{\T} n\,\p_t\varphi\,\diff x\diff t
 -\int_{\T}n_0(x)\varphi(0,x)\,\diff x
 +\nu\int_0^T  \int_{\T} n\partial_x\rho\,\p_x\varphi\,\diff x\diff t=0.
\end{equation}
\end{enumerate}
\end{definition}
\begin{thm}
\label{thm:main}
Let $T>0$, $\nu>0$, $\nu\neq 1$, and initial data $m_0,n_0\in L^\infty(\T)$, $m_0,n_0\ge 0$. Let $(m^k,n^k)$ be an admissible approximation sequence in the sense of Definition~\ref{def:admissible}. Then, up to a subsequence not relabeled, $(m^k,n^k)$ converges to $(m,n)$ which is a weak solution of \eqref{eq:MNsystem} on $(0,T)\times\T$ in the sense of Definition~\ref{def:weak-solution}.
\end{thm}

\begin{remark}
The viscous approximation, where one adds $\eps\partial_{xx}m,\,  \eps\partial_{xx}n$ and pass to the limit $\eps\to 0$ yields an admissible approximation sequence.
\end{remark}

The proof has three main steps. First, we build an entropy family. This family is then used together with the div--curl lemma to control the correlations between $a^k$ and $\partial_x\rho^k$. Then a second commutation identity shows that the Young measure associated to $a$ is a Dirac mass wherever the pressure gradient is nonzero. The extension of this result to the higher dimensional case, or more than two species is left as open.

\section{\texorpdfstring{Structure of $\rho$ and $a$}{Structure of rho and a}}
\label{sec:structure}

The goal of this section is to obtain estimates on $\rho,a$ and derive the balance laws necessary to obtain enough compactness. In a first subsection, we focus on basic estimates, taking advantage of the change of variable and the structure of the equations. Then, after introducing some notations, we show that the equations admit in fact an infinite number of balance laws. 

\subsection{Evolution and first estimates}

We introduce the nonnegative polynomial
\begin{equation}
\label{eq:defB}
 B(a):=(a-1)(\nu-a)=(a-\al)(\be-a),
 \quad a\in I.
\end{equation}
Also we have the identity $\al\be=\nu$. The following proposition follows from a simple computation, and sets up the basic identities on the variables $(\rho,a)$.

\begin{proposition}
\label{prop:rhoa-system}
Let $(m,n)$ be a smooth positive solution of \eqref{eq:MNsystem}, and let $(\rho,a)$ be defined by \eqref{eq:mnfromrhoa}. Then the total density satisfies
\begin{equation}
\label{eq:rhoeq}
 \p_t\rho - \p_x\left(a\rho\p_x\rho\right)=0.
\end{equation}
The product of the density and the activity variable satisfies
\begin{equation}\label{eq:arhoeq}
\partial_t (a\rho ) - \partial_x(((\nu+1) a -\nu)\rho \partial_x\rho)=0.    
\end{equation}
Moreover, the activity variable solves
\begin{equation}
\label{eq:aeq}
 \p_t a - (\nu+1-a) \partial_x\rho \partial_x a
 = B(a)\left(\partial_{xx}\rho+\frac{(\partial_x\rho)^2}{\rho}\right).
\end{equation}

\end{proposition}

Equation~\eqref{eq:aeq} is not really useful by itself, and we only use it later implicitly to derive the family of balance laws.

\begin{proof}
Summing the two equations in \eqref{eq:MNsystem} gives
$$
 \p_t\rho - \p_x\left((m+\nu n) \partial_x\rho\right)=0.
$$
Since $m+\nu n=a\rho$, we obtain \eqref{eq:rhoeq}. Then we use the identity  $m+\nu^2n=((\nu+1)a-\nu)\rho$ and we obtain~\eqref{eq:arhoeq} in a similar way. To derive the equation for $a$, we let $\omega=m+\frac1\nu n.$ Then
$$
 \partial_t\omega - \p_x(\rho \partial_x\rho)=0, \quad
 \omega=q\rho, \quad q=\frac{\omega}{\rho}=\frac{\nu+1-a}{\nu}.
$$
Using also \eqref{eq:rhoeq}, we find
$$
 \partial_t q-\nu q \partial_x\rho \partial_x q=(1-q)(1-\nu q)\left(\partial_{xx}\rho+\frac{ (\partial_x\rho)^2}{\rho}\right).
$$
Multiplying by $-\nu$ and using $a=\nu+1-\nu q$, $ \partial_x a=-\nu \partial_x q$, and
$$
 -\nu(1-q)(1-\nu q)=(a-1)(\nu-a)=B(a),
$$
we obtain \eqref{eq:aeq}. 

\end{proof}

The following proposition follows from a simple computation. It yields the basic estimates used later for weak compactness of $m^k$, $n^k$, and $\partial_x\rho^k$, together with strong compactness of $\rho^k$.

\begin{proposition}
\label{prop:basic-identities}
Let $(m,n)$ be a smooth positive solution of \eqref{eq:MNsystem}. Then:
\begin{enumerate}[label=(\roman*),leftmargin=2.4em]
\item We have the entropy identity
\begin{equation}
\label{eq:basic-entropy}
\frac{\diff}{\diff t}
\int_{\T}\left(h(m)+\frac1\nu h(n)\right) \diff x
+
\int_{\T}| \partial_x\rho|^2 \diff x=0,
\end{equation}
where $h(z)=z\log z-z$.
\item The maximum principle holds, i.e.
\begin{equation}
\label{eq:maxprinciple-rho}
 \|\rho(t,\cdot)\|_{L^\infty(\T)}\leq \|\rho(0,\cdot)\|_{L^\infty(\T)}
 \quad \text{for all } t\ge 0.
\end{equation}
\end{enumerate}
Moreover these estimates hold uniformly for a sequence of smooth solutions $(m^k, n^k)$ of System~\eqref{eq:MNsystem} with $L^{\infty}(\T)$ initial conditions.
\end{proposition}

The maximum principle can be seen for instance after multiplying the equation~\eqref{eq:rhoeq} by $\rho^p$, integrating in space and then sending $p\to +\infty$. These two estimates control the total density.

\subsection{The balance laws family}

In order to find some new balance laws for this system, we first introduce some notations.

\begin{lemma}
\label{lem:X-transform}
Let $X:I^\circ\to\R$ be defined by
\begin{equation}
\label{eq:defX}
 X'(a)=\frac{a}{B(a)}=\frac{a}{(a-\al)(\be-a)}.
\end{equation}
Then $X$ is an increasing diffeomorphism from $I^\circ=(\al,\be)$ onto $\R$. Moreover,
\begin{equation*}
 X(a)=\frac{\al}{\be-\al}\log(a-\al)-\frac{\be}{\be-\al}\log(\be-a)+C,
\end{equation*}
for a suitable constant $C$.
\end{lemma}

\begin{lemma}
\label{lem:phi-family}
Let the interval 
\begin{equation*}
 S:=\left(\frac{\al}{\be},\frac{\be}{\al}\right)
 =\left(\min\{\nu,\nu^{-1}\},\max\{\nu,\nu^{-1}\}\right).
\end{equation*}
For every $s\in S$, define $\phi_s$ on $I$ by
\begin{equation*}
 \phi_s(a):=\int_{\al}^{a} e^{-(s-1)X(r)} \diff r,
 \quad a\in I.
\end{equation*}
Then $\phi_s\in C(I)\cap C^2(I^\circ)$, $\phi_s(\al)=0$, and
\begin{equation}
\label{eq:phi-ode}
 B(a)\phi_s''(a)+(s-1)a\phi_s'(a)=0
 \quad \text{for } a\in I^\circ.
\end{equation}
Moreover,
\begin{equation}
\label{eq:Ms-def}
 M_s(a):=s a\phi_s(a)+B(a)\phi_s'(a)
\end{equation}
is continuous on $I$, belongs to $C^1(I^\circ)$, and satisfies
\begin{equation}
\label{eq:Ma-derivative}
M_s'(a)=s\phi_s(a)+(\nu+1-a)\phi_s'(a), \quad \left(\frac{M_s(a)}{a}\right)'=\frac{\nu}{a^2}\phi_s'(a)
 \quad \text{for } a\in I^\circ.
\end{equation}
\end{lemma}

\begin{proof}
By construction,
$$
\phi_s'(a)=e^{-(s-1)X(a)},
\quad
\phi_s''(a)=-(s-1)X'(a)e^{-(s-1)X(a)}.
$$
Using \eqref{eq:defX}, we obtain \eqref{eq:phi-ode}. The continuity of $\phi_s$ at the endpoints follows from the definition of $S$. Indeed, by Lemma~\ref{lem:X-transform},
$$
\phi_s'(a)\sim (a-\al)^{-\frac{(s-1)\al}{\be-\al}}
\quad \text{as } a\to \al,
\quad
\phi_s'(a)\sim (\be-a)^{\frac{(s-1)\be}{\be-\al}}
\quad \text{as } a\to\be,
$$
and both singularities are integrable exactly when $s\in (\al/\be,\be/\al)$.

We now prove \eqref{eq:Ma-derivative}. We differentiate \eqref{eq:Ms-def} and use the ODE~\eqref{eq:phi-ode}, as well as the identity $B'(a)=\nu+1-2a$:
\begin{align*}
 M_s'(a)
 &= s\phi_s(a)+sa\phi_s'(a)+B'(a)\phi_s'(a)+B(a)\phi_s''(a)\\
 &= s\phi_s(a)+\left(sa+B'(a)-(s-1)a\right)\phi_s'(a)\\
 &= s\phi_s(a)+(\nu+1-a)\phi_s'(a),
\end{align*}
This yields
\begin{align*}
\left(\frac{M_s}{a}\right)'
&=\frac{a M_s'-M_s}{a^2}\\
&=\frac{a\left(s\phi_s+(\nu+1-a)\phi_s'\right)-\left(sa\phi_s+B\phi_s'\right)}{a^2}\\
&=\frac{a(\nu+1-a)-B(a)}{a^2} \phi_s'(a).
\end{align*}
Since $B(a)=-(a^2-(\nu+1)a+\nu)$, the numerator equals $\nu$, and \eqref{eq:Ma-derivative} follows.
\end{proof}

With the previous notations and identities, we introduce the new family of balance laws, that we use later to prove strong compactness with the div--curl lemma. The functions $\phi_s$, $M_s$ are designed so that the bad $\partial_x a$ and $\partial_{xx}\rho$ terms cancel in the next computation. Note that we don't dissipate $\rho\phi_s(a)$ and we rather dissipate $\rho^s\phi_s(a)$ since the first choice would lead to the condition $\phi_s''(a)=0$ and therefore $\phi_s$ is affine, an entropy we already know.

\begin{proposition}
\label{prop:entropy-identity}
Let $(m,n)$ be a smooth positive solution of \eqref{eq:MNsystem}, let $(\rho,a)$ be respectively the sum and activity variable of $(m,n)$, and let $s\in S$. Then
\begin{equation}
\label{eq:exact-balance-law}
 \p_t\left(\rho^s\phi_s(a)\right)
 - \p_x\left(M_s(a)\rho^s \partial_x\rho\right)
 =(1-s)M_s(a)\rho^{s-1} (\partial_x\rho)^2.
\end{equation}
\end{proposition}

\begin{proof}
We define
$$
U_s=\rho^s\phi_s(a),
\quad
F_s=M_s(a)\rho^s \partial_x\rho.
$$
By differentiating in time we obtain:
\begin{align*}
\p_t U_s
&= s\rho^{s-1}\phi_s(a)\partial_t \rho+\rho^s\phi_s'(a)\partial_t a.
\end{align*}
Using \eqref{eq:rhoeq} and \eqref{eq:aeq}, and with $\xi= \partial_x\rho$, we obtain
\begin{align*}
\p_t U_s
&= s\rho^{s-1}\phi_s \p_x(a\rho\xi)
 +\rho^s\phi_s'\left((\nu+1-a)\xi  \partial_x a + B(a)\left(\partial_x\xi+\xi^2/\rho\right)\right)\\
&=\left(sa\phi_s+B\phi_s'\right)\rho^s\partial_x\xi
 +\left(sa\phi_s+B\phi_s'\right)\rho^{s-1}\xi^2
 +\rho^s\xi  \partial_x a\left(s\phi_s+(\nu+1-a)\phi_s'\right).
\end{align*}
On the other hand,
\begin{align*}
\p_x F_s
&=\rho^s\xi  \partial_x a M_s'(a)+sM_s(a)\rho^{s-1}\xi^2+M_s(a)\rho^s\partial_x\xi.
\end{align*}
Subtracting, and using the formula for $M_s'$ from Lemma~\ref{lem:phi-family}, the $ \partial_x a$-terms cancel. Hence
\begin{align*}
\p_t U_s-\p_x F_s
&=\left(M_s\rho^{s-1}\xi^2-sM_s\rho^{s-1}\xi^2\right)\\
&=(1-s)M_s(a)\rho^{s-1} (\partial_x\rho)^2.
\end{align*}
This proves \eqref{eq:exact-balance-law}.
\end{proof}

We show that the family $\{\phi_s\}$ is dense into $C(I)$. This will allow us to extend the compensated-compactness identities from the entropy family to all continuous functions of $a$.

\begin{proposition}
\label{prop:density}
Let $J\subset S$ be a nonempty open interval. Then the linear span of
$$
\{1\}\cup\{\phi_s:s\in J\}
$$
is dense in $C(I)$ with respect to the uniform norm.
\end{proposition}

\begin{figure}[H]
    \centering
    \includegraphics[scale=0.5]{Figure_2.png}
    \caption{Sample plots of $\phi_s$ for $s\in\{0.7,0.8,\dots,1.3\}$, with $\nu=2$, $I=[1,2]$, and $S=(\frac12,2)$.}
    \label{fig}
\end{figure}

\begin{proof}
It is enough to show that if a finite signed measure $\eta$ on $I$ satisfies
$$
 \int_I 1 \diff\eta=0,
 \quad
 \int_I \phi_s(a) \diff\eta(a)=0
 \quad \forall s\in J,
$$
then $\eta=0$.

We define the tail function 
$$
 G(r)=\eta([r,\be]),
 \quad r\in I.
$$
The goal is to prove that $G$ vanishes for all $r\in I$. This yields $\eta=0$. 
Since $\eta(I)=0$, we already have $G(\al)=0$. By Fubini's theorem and the definition of $\phi_s$,
\begin{align*}
 0 &=\int_I \phi_s(a) \diff\eta(a)  =\int_I\int_{\al}^{a} e^{-(s-1)X(r)} \diff r \diff\eta(a)\\
 &=\int_{\al}^{\be} e^{-(s-1)X(r)} G(r) \diff r,
 \quad \forall s\in J.
\end{align*}
Now let $x=X(r)$, so that $r=X^{-1}(x)$ and
$$
 h(x)=G(X^{-1}(x))(X^{-1})'(x).
$$
Since $X:I^\circ\to\R$ is a diffeomorphism and $J\subset S$, the function $e^{-(s-1)x}h(x)$ is integrable for every $s\in J$. Therefore the Laplace transform
$$
 F(s)=\int_{\R} e^{-(s-1)x}h(x) \diff x
$$
is holomorphic in the vertical strip $\{z\in\C: \Re z\in J\}$, and vanishes on the real interval $J$. Hence $F\equiv 0$ on that strip. Now let $s_0\in J$. For every $\tau\in\R$,
$$
 0=F(s_0+i\tau)
 =\int_{\R} e^{-i\tau x}e^{-(s_0-1)x}h(x) \diff x.
$$
Thus the Fourier transform of the $L^1$-function $e^{-(s_0-1)x}h(x)$ vanishes identically. Hence $h=0$ almost everywhere, and therefore $G=0$ almost everywhere on $I$.

The function $G$ is right-continuous. Since a right-continuous function which vanishes almost everywhere is identically zero, we conclude that $G\equiv 0$. Finally,  $\eta=0$.
\end{proof}

\section{Admissible approximations and a first compactness reduction}
\label{sec:young}

We now turn from smooth solutions to approximation sequences. At this point the total density is compact, but the activity variable may still oscillate. The aim of this section is to show that those oscillations already satisfy a particular constraint.

\subsection{Admissible sequences, div--curl lemma and Young measures}

\begin{definition}
\label{def:admissible}
Let $T>0$, $\nu>0$, $\nu\neq 1$, and choose a nonempty open interval $J\subset(1,\be/\al)$. A nonnegative sequence $(m^k,n^k)_{k\in\N}$ on $[0,T)\times\T$ is called an \emph{admissible approximation sequence} if the following properties hold.
\begin{enumerate}[label=(A\arabic*)]
\item \textit{Positivity and boundedness.}\label{A1}
Define
$$
 \rho^k=m^k+n^k,
 \quad
 a^k=\frac{m^k+\nu n^k}{\rho^k}.
$$
Then $\rho^k> 0$, $\rho^k$ is bounded uniformly in $L^\infty((0,T)\times\T)$, and $a^k(t,x)\in I=[\al,\be]$ for almost every $(t,x)$.

\item \textit{Compactness of the total density.}\label{A2}
There exists a nonnegative function $\rho$ such that for all $1\le p<+\infty$
$$
 \rho^k\to \rho \quad \text{strongly in } L^p((0,T)\times\T),
$$
and
$$
 \xi^k=\p_x\rho^k \weakto \p_x\rho \quad \text{weakly in } L^2((0,T)\times\T).
$$

\item \textit{Weak compactness of the species.}\label{A3}
There exist nonnegative functions $m,n\in L^\infty((0,T)\times\T)$ such that, 
$$
 m^k \overset{*}{\rightharpoonup} m,
 \quad
 n^k \overset{*}{\rightharpoonup} n, \quad a^k \overset{*}{\rightharpoonup} a
 \quad \text{weakly star in } L^\infty((0,T)\times\T),
$$

Moreover $m+n=\rho$ almost everywhere.

\item \textit{Affine balance laws and consistency.}\label{A4}
The two sequences of vector fields
$$
X_0^k=\left(\rho^k,-a^k\rho^k\xi^k\right),
\quad
X_1^k=\left(a^k\rho^k,-((\nu+1)a^k-\nu)\rho^k\xi^k\right)
$$
are bounded in $L^2((0,T)\times\T;\R^2)$. Their residuals
$$
r_0^k=\p_t\rho^k-\p_x(a^k\rho^k\xi^k),
\quad
r_1^k=\p_t(a^k\rho^k)-\p_x\left(((\nu+1)a^k-\nu)\rho^k\xi^k\right)
$$
satisfy
$$
r_0^k\to 0,
\quad
r_1^k\to 0
\quad \text{in } L^2(0,T;H^{-1}(\T)).
$$

\item \textit{Family of balance laws.}\label{A5}
For every $s\in J$, the vector field
$$
X_s^k=\left((\rho^k)^s\phi_s(a^k),-M_s(a^k)(\rho^k)^s\xi^k\right)
$$
is bounded in $L^2((0,T)\times\T;\R^2)$, and its residual
$$
r_s^k=\p_t\left((\rho^k)^s\phi_s(a^k)\right)-\p_x\left(M_s(a^k)(\rho^k)^s\xi^k\right)
$$
is bounded in $L^1((0,T)\times\T)$.

\item \textit{Initial data.}\label{A6}
There exist $m_0,n_0\in L^\infty(\T)$, $m_0,n_0\ge 0$, such that
$$
 m^k(0,\cdot)\to m_0,
 \quad
 n^k(0,\cdot)\to n_0
 \quad \text{strongly in } L^2(\T).
$$
\end{enumerate}
\end{definition}

These assumptions are designed for approximation schemes which preserve the balance laws derived above for smooth solutions. Indeed, the previous computations show that a sequence of smooth solutions satisfy exactly these estimates. For instance, \ref{A1}--\ref{A2}--\ref{A3} follow from Proposition~\ref{prop:basic-identities} and Aubin-Lions Lemma. \ref{A4} follows from $r_0^k = r_1^k=0$. Finally \ref{A5} follows from continuity of the functions $M_s, \phi_s$ on the compact interval $I$.
In other words, the admissibility assumptions simply record the estimates and balance laws already present in the equation. To construct an approximation scheme, one can for instance add an artificial viscous regularization with a parameter $\eps>0$ and then pass to the inviscid limit $\eps\to 0$. Let us also mention that the choice to use an admissible sequence and not design a particular scheme is motivated by the fact that we do not know whether uniqueness of System~\ref{eq:MNsystem} should be expected. Therefore the constructed solutions would be associated to a particular scheme. Our result just says that if this scheme preserves the structure of the equation, it yields existence of a solution.

The new family of balance laws by itself do not provide any compactness estimates. One needs to use the div--curl lemma, now a common tool in conservation laws. In our case, we use a consequence of the div--curl lemma. More precisely, we use the version from~\cite{Skrzeczkowski2026}, which can be deduced from~\cite{MR2769903}:

\begin{lemma}[Consequence of the div--curl lemma]
\label{lem:tartar}
Let $(U^k,Q^k)$ and $(V^k,P^k)$ be such that 
\begin{enumerate}
    \item $(U^k)$ is bounded uniformly in $L^\infty(0,T; L^{2}(\T))$
    \item $(V^k), (Q^k), (P^k)$ are bounded uniformly in $L^2(0,T; L^{2}(\T))$
    \item $\p_t U^k+\p_x Q^k$ are precompact in $L^2(0,T; H^{-1}(\T))$
    \item $\p_t V^k+\p_x P^k$ is bounded in $L^1(0,T; L^{1}(\T))$
\end{enumerate}
Then, after passage to a subsequence,
\begin{equation}
\label{eq:tartar-identity}
 \overline{U P-Q V}=\overline U\, \overline P-\overline Q\, \overline V
 \quad \text{in } \mathcal D'((0,T)\times\T),
\end{equation}
where the bar denotes weak limits.
\end{lemma}

Let $(m^k,n^k)$ be an admissible approximation sequence in the sense of Definition~\ref{def:admissible}. By boundedness of $a^k$ in $L^\infty$ and of $\xi^k$ in $L^2$, there exists, after passage to a subsequence, a parametrized Young measure
$$
 \la_{t,x}\in \mathcal P(I\times\R)
 \quad \text{for a.e. } (t,x)\in (0,T)\times\T,
$$
generated by the pair $(a^k,\xi^k)$. We denote by $\mu_{t,x}$ the marginal of $\la_{t,x}$ on $I$:
$$
 \mu_{t,x}(E):=\la_{t,x}(E\times \R),
 \quad E\subset I\ \text{Borel}.
$$
After redefining on a null set, we may also write the disintegration
$$
 \la_{t,x}(\diff a,\diff\xi)=\mu_{t,x}(\diff a)\,\la^a_{t,x}(\diff\xi).
$$
Whenever no confusion is possible, we simply write $\av{\cdot}$ for integration with respect to $\la_{t,x}$ and $\avmu{\cdot}$ for integration with respect to $\mu_{t,x}$. With all the ingredients let us explain the main strategy, that is also sketched in the introduction. In order to pass to the limit we need to pass to the limit in the nonlinear flux
$$
m^k\partial_x\rho^k = \frac{\nu}{\nu-1}\rho^k\partial_x\rho^k - \frac{1}{\nu-1}a^k\rho^k\partial_x\rho^k.
$$

As $\rho^k$ converges strongly we see that the only difficulty is to prove that $a^k\partial_x\rho^k$ converges weakly to $a\partial_x\rho$. Or, in the Young measure language, that $\langle a \xi\rangle = a\partial_x\rho$. Putting aside the technical difficulties generated by the fact that $\rho$ is only nonnegative and not positive, proving that $\mu_{t,x}$ is a Dirac mass yields 
$$
\langle a\xi\rangle = a\langle\xi\rangle = a\partial_x\rho.
$$
The rest of the paper is thus devoted to proving that $\mu_{t,x}$ is a Dirac mass.
The following lemma is a simple consequence of the strong convergence of $\rho^k$ and weak-* convergence of $a^k$.

\begin{lemma}
\label{lem:weak-limits-via-young}
For every bounded continuous function $f:I\to\R$, up to a subsequence,
$$
 f(a^k)\rho^k \weakto \rho \avmu{f(a)}
 \quad \text{weakly in } L^2((0,T)\times\T).
$$
Also, since $a^k$ converges weakly star and $\rho^k$ converges strongly:

\begin{equation*}
 a^k\rho^k \weakto \rho a
 \quad \text{weakly in } L^2((0,T)\times\T).
\end{equation*}
Moreover,
\begin{equation*}
 \rho a = m+\nu n
 \quad \text{a.e. in } (0,T)\times\T.
\end{equation*}
\end{lemma}

\subsection{The first div--curl identity}
The next proposition shows that we can find a first structure between the oscillations of $a^k$ and $\xi^k$.

\begin{proposition}
\label{prop:first-hit}
For almost every $(t,x)\in(0,T)\times\T$ such that $\rho(t,x)>0$, and for every $s\in J$, let
$$
A=A(t,x):=\av{a\xi}.
$$
Then the following properties hold:
\begin{enumerate}
    \item There exists a countable dense subset $D\subset J$, such that for every $s\in D$ 
    \begin{equation}
    \label{eq:first-hit-entropy}
        \av{a\xi \phi_s(a)}
        = A \avmu{\phi_s(a)}.
    \end{equation}

    \item Since $\{1\}\cup\{\phi_s:s\in J\}$ is dense in $C(I)$, and $s\mapsto \phi_s$ is continuous as a map from $J\to C(I)$, for every $f\in C(I)$,
    \begin{equation}
    \label{eq:first-hit-allf}
        \av{a\xi f(a)}
        = A \avmu{f(a)}.
    \end{equation}

    \item For $\mu_{t,x}$-a.e. $a$,
    \begin{equation}
    \label{eq:conditional-xi}
        \int_{\R}\xi\,\lambda^a_{t,x}(\diff\xi)=\frac{A(t,x)}{a}.
    \end{equation}

    \item One has
    \begin{equation}
    \label{eq:rhox-A}
        \partial_x\rho = A \avmu{\frac1a}.
    \end{equation}

    \item On the set $\{\rho>0\}$,
    \begin{equation}
    \label{eq:Aiff}
        A=0 \quad \Longleftrightarrow \quad \partial_x\rho=0.
    \end{equation}
\end{enumerate}
\end{proposition}

\begin{proof}
Let $s\in J$. We apply Lemma~\ref{lem:tartar} to the two sequences
$$
 (U^k,Q^k)=\left(\rho^k,-a^k\rho^k\xi^k\right),
 \quad
 (V^k,P^k)=\left((\rho^k)^s\phi_s(a^k),-M_s(a^k)(\rho^k)^s\xi^k\right).
$$
By admissibility we obtain
$$
 \overline{U P-Q V}=\overline U\, \overline P-\overline Q\, \overline V.
$$
Using the strong convergence of $\rho^k$ and the Young measure $\la$, we obtain 
\begin{align*}
 \overline{U P-Q V}
 &=\rho^{s+1}\av{(a\phi_s(a)-M_s(a))\xi},\\
 \overline U\, \overline P-\overline Q\, \overline V
 &=\rho^{s+1}\left(A \avmu{\phi_s(a)}-\av{M_s(a)\xi}\right).
\end{align*}
Therefore,
$$
 \rho^{s+1}\left(\av{a\xi \phi_s(a)}-A \avmu{\phi_s(a)}\right)=0.
$$
Since $\rho>0$ on the set under consideration, \eqref{eq:first-hit-entropy} holds a.e. for this particular $s$. This can then be extended to a countable subset $D$.

Since the span of $\{1\}\cup\{\phi_s:s\in J\}$ is dense in $C(I)$ and $s\mapsto \phi_s$ is continuous as a map from $J\to C(I)$, the functional
$$
 f\mapsto \av{a\xi f(a)}-A \avmu{f(a)}
$$
which is continuous vanishes on a dense subspace of $C(I)$, hence on all of $C(I)$. This proves \eqref{eq:first-hit-allf}. Finally, \eqref{eq:first-hit-allf} is equivalent to
$$
 \int_I f(a)\left(a\int_{\R}\xi \la^a_{t,x}(\diff\xi)-A\right)\mu_{t,x}(\diff a)=0
 \quad \forall f\in C(I).
$$
Thus the bracket vanishes for $\mu_{t,x}$-almost every $a$. Dividing by $a\in I\subset (0,\infty)$ yields \eqref{eq:conditional-xi}. Integrating \eqref{eq:conditional-xi} over $\mu_{t,x}$ gives \eqref{eq:rhox-A}. Then, \eqref{eq:Aiff} follows.
\end{proof}

From~\eqref{eq:Aiff} and the assumptions of Proposition~\ref{prop:first-hit}, we see that the sets $\{\rho>0\}$ and $\{\partial_x\rho\ne 0\}$ play an important role.  The following proposition will be therefore useful and follows from~\cite[Chapter 6]{MR1158660}, adapted to space-time functions.

\begin{lemma}\label{lem:vacuum-gradient}
Assume $\rho\in L^2(0,T;H^1(\T))$ and $\rho\ge 0$ a.e. Then
$$
\{\partial_x\rho\neq 0\}\subset \{\rho>0\}
\quad\text{up to a null set.}
$$
\end{lemma}

\section{Collapse of the Young measure}

It remains to show that the measure $\mu_{t,x}$ is a Dirac mass on the set where the pressure gradient is not zero. This comes from a second compensated-compactness identity which acts only on $\mu$ of the variable $a$.

\subsection{The second identity}

\begin{proposition}
\label{prop:commutator-on-entropy-family}
Let $(t,x)$ be a point such that $\rho(t,x)>0$ and $A(t,x)\neq 0$, where $A=\av{a\xi}$. Then there exists a countable dense subset $D\subset J$, such that for every $s\in D$,
\begin{equation}
\label{eq:entropy-commutator}
\Cov_\mu \left(a,\frac{M_s(a)}{a}\right)
=
\Cov_\mu \left(\phi_s(a),-\frac{\nu}{a}\right).
\end{equation}
\end{proposition}

\begin{proof}
We apply Lemma~\ref{lem:tartar} to the affine balance law for $a\rho$ and to the entropy balance indexed by $s$:
$$
 (U^k,Q^k)=\left(a^k\rho^k,-((\nu+1)a^k-\nu)\rho^k\xi^k\right),
$$
$$
 (V^k,P^k)=\left((\rho^k)^s\phi_s(a^k),-M_s(a^k)(\rho^k)^s\xi^k\right).
$$
By admissibility, the divergence of $(U^k,Q^k)$ converges to $0$ in $L^2(0,T;H^{-1}(\T))$, while the divergence of $(V^k,P^k)$ is bounded in $L^1((0,T)\times\T)$. Hence Lemma~\ref{lem:tartar} applies. Using the strong convergence of $\rho^k$, we obtain
\begin{align*}
 \overline{U P-Q V}
 &=\rho^{s+1}\left(
 -\av{aM_s(a)\xi}
 +\av{((\nu+1)a-\nu)\phi_s(a)\xi}
 \right),\\
 \overline U\, \overline P-\overline Q\, \overline V
 &=\rho^{s+1}\left(
 -a\,\av{M_s(a)\xi}
 +\bigl((\nu+1)A-\nu \partial_x\rho\bigr)\avmu{\phi_s(a)}
 \right).
\end{align*}
Since $\rho>0$, with Proposition~\ref{prop:first-hit} we obtain
\begin{align*}
 \av{aM_s(a)\xi}&=A\avmu{M_s(a)},\\
 \av{((\nu+1)a-\nu)\phi_s(a)\xi}
 &=A\left((\nu+1)\avmu{\phi_s(a)}-\nu\avmu{\frac{\phi_s(a)}{a}}\right),\\
 \av{M_s(a)\xi}&=A\avmu{\frac{M_s(a)}{a}},\\
 \partial_x\rho&=A\avmu{\frac1a}.
\end{align*}
Substituting these identities into the div--curl relation, dividing by $\rho^{s+1}A$, and simplifying, we obtain
\begin{equation*}
 \avmu{M_s(a)}-a\,\avmu{\frac{M_s(a)}{a}}
 =
 -\nu\left(
 \avmu{\frac{\phi_s(a)}{a}}
 -
 \avmu{\phi_s(a)}\avmu{\frac1a}
 \right).
\end{equation*}
The previous identity can be rewritten as
$$
 \Cov_\mu\left(a,\frac{M_s}{a}\right)
 =
 -\nu\,\Cov_\mu\left(\phi_s,\frac1a\right).
$$
Thus, we obtain \eqref{eq:entropy-commutator}.
\end{proof}

\subsection{Collapse on the nondegenerate set}
\label{sec:collapse}
In this section we collapse the scalar Young measure $\mu$ to a Dirac mass on the nondegenerate set $A(t,x)\neq0$, and $\rho(t,x)>0$.

\begin{thm}
\label{thm:collapse}
For almost every $(t,x)\in(0,T)\times\T$ such that $A(t,x)\neq0$, and $\rho(t,x)>0$ the measure $\mu_{t,x}$ is a Dirac mass at $a(t,x)$.
\end{thm}

\begin{proof}
Let $(t,x)$ such that $A(t,x)\neq0$ and $\rho(t,x)>0$, and write $\mu=\mu_{t,x}$.
By Proposition~\ref{prop:commutator-on-entropy-family}, for some $s\in J$,
$$
\Cov_\mu \left(a,\frac{M_s(a)}{a}\right)
=
\Cov_\mu \left(\phi_s(a),-\frac{\nu}{a}\right).
$$
We recall the following identity for the covariance: 
$$
\Cov_\mu(f,h)
=
\frac12\int_I\int_I (f(r)-f(q))(h(r)-h(q))\,d\mu(r)\,d\mu(q).
$$
Applying the previous identity and symmetry we obtain
$$
0=
\iint_{\{r<q\}}
\left(
(q-r)\left(\frac{M_s(q)}{q}-\frac{M_s(r)}{r}\right)
-\nu(\phi_s(q)-\phi_s(r))\left(\frac1r-\frac1q\right)
\right)\,d\mu(r)\,d\mu(q).
$$
By \eqref{eq:Ma-derivative},
$$
\frac{M_s(q)}{q}-\frac{M_s(r)}{r}
=
\nu\int_r^q \frac{\phi_s'(u)}{u^2}\,du.
$$
Also
$$
\phi_s(q)-\phi_s(r)=\int_r^q \phi_s'(u)\,du,
\qquad
\frac1r-\frac1q=\int_r^q \frac{1}{u^2}\,du.
$$
Therefore
$$
0=
\nu\iint_{\{r<q\}}
\left(
(q-r)\int_r^q \frac{\phi_s'(u)}{u^2}\,du
-
\left(\int_r^q \phi_s'(u)\,du\right)
\left(\int_r^q \frac{1}{u^2}\,du\right)
\right)\,d\mu(r)\,d\mu(q).
$$
Since $s\in J\subset (1,\beta/\alpha)$ and $X$ is increasing, the function
$$
u\mapsto \phi_s'(u)=e^{-(s-1)X(u)}
$$
is strictly decreasing on $I$. Also the function $u\mapsto 1/u^2$ is strictly decreasing on $I$. Hence, by Chebyshev's integral inequality,
$$
(q-r)\int_r^q \frac{\phi_s'(u)}{u^2}\,du
>
\left(\int_r^q \phi_s'(u)\,du\right)
\left(\int_r^q \frac{1}{u^2}\,du\right)
\quad \text{for every } r<q.
$$
Thus the integrand is strictly positive on $\{r<q\}$.

If $\mu$ were not a Dirac mass, then $\mu\otimes\mu(\{r<q\})>0$, so the integral would be strictly positive, which is a contradiction. Hence $\mu$ is a Dirac mass. By weak convergence of $a^k$ we identify this Dirac as $\delta_{a(t,x)}$.
\end{proof}
\begin{corollary}
\label{cor:A-equals-abar-rhox}
We have
\begin{equation}
\label{eq:A-equals-abar-rhox}
 \rho A = \rho a  \partial_x\rho
 \quad \text{a.e. in } (0,T)\times\T.
\end{equation}
Therefore,
\begin{equation}
\label{eq:limit-j}
 a^k\rho^k\xi^k \weakto \rho a  \partial_x\rho
 \quad \text{weakly in } L^1((0,T)\times\T).
\end{equation}
\end{corollary}

\begin{proof}
On the set $\{ \partial_x\rho\neq 0\}$, Lemma~\ref{lem:vacuum-gradient} gives $\rho>0$, Equation~\eqref{eq:Aiff} gives $A\ne 0$ and Theorem~\ref{thm:collapse} shows that $\mu$ is a Dirac mass at $a$. Hence
$$
 A=\av{a\xi}=a \av{\xi}=a  \partial_x\rho.
$$
On the set $\{\rho>0\}\cap\{\partial_x\rho=0\}$, Proposition~\ref{prop:first-hit} gives $A=0$. Multiplying by $\rho$ proves \eqref{eq:A-equals-abar-rhox}; on the set $\{\rho=0\}$ both sides vanish. The weak convergence \eqref{eq:limit-j} follows from the Young-measure representation of the weak limit of $a^k\rho^k\xi^k$:
$$
 a^k\rho^k\xi^k \weakto \rho \av{a\xi}=\rho A=\rho a  \partial_x\rho.
$$
\end{proof}

\section{Identification of the fluxes}
\label{sec:existence}

We can now identify the limits of the two nonlinear fluxes.

\begin{lemma}
\label{lem:flux-identification}
Define
$$
 J_m^k:=m^k\xi^k,
 \quad
 J_n^k:=n^k\xi^k.
$$
Then, up to a subsequence not relabeled,
\begin{equation}
\label{eq:flux-m-limit}
 J_m^k \weakto  m  \partial_x\rho,
 \quad   J_n^k \weakto  n  \partial_x\rho\quad \text{weakly in } L^1((0,T)\times\T),
\end{equation}
\end{lemma}

\begin{proof}
We recall that the two species can be rewritten in terms of the total density and the activity variable
$$
 m^k=\frac{\nu-a^k}{\nu-1}\rho^k,
 \quad
 n^k=\frac{a^k-1}{\nu-1}\rho^k.
$$
Thus, we can rewrite the flux as 
\begin{equation}
\label{eq:flux-decomp-m}
 J_m^k
 =\frac{\nu\rho^k\xi^k-a^k\rho^k\xi^k}{\nu-1}, \quad  J_n^k
 =\frac{a^k\rho^k\xi^k-\rho^k\xi^k}{\nu-1}.
\end{equation}

Now $\rho^k\to\rho$ strongly in $L^2$ and $\xi^k\weakto  \partial_x\rho$ weakly in $L^2$, hence
$$
 \rho^k\xi^k \weakto \rho \partial_x\rho
 \quad \text{weakly in } L^1((0,T)\times\T).
$$
Moreover, by Corollary~\ref{cor:A-equals-abar-rhox},
$$
 a^k\rho^k\xi^k \weakto \rho a  \partial_x\rho
 \quad \text{weakly in } L^1((0,T)\times\T).
$$
We can therefore identify the flux: 
$$
 J_m^k \weakto \frac{\nu\rho \partial_x\rho-\rho a  \partial_x\rho}{\nu-1},
 \quad
 J_n^k \weakto \frac{\rho a  \partial_x\rho-\rho \partial_x\rho}{\nu-1}.
$$
Since $\rho= m+n$ and $\rho a= m+\nu n$ by Lemma~\ref{lem:weak-limits-via-young}, we obtain the result.
\end{proof}

We can now turn to the proof of the main theorem. 

\begin{proof}[Proof of Theorem~\ref{thm:main}]
Let
$$
r_0^k=\p_t\rho^k-\p_x(a^k\rho^k\xi^k),
\quad
r_1^k=\p_t(a^k\rho^k)-\p_x\left(((\nu+1)a^k-\nu)\rho^k\xi^k\right).
$$
By Assumption~\ref{A4}, both residuals converge to $0$ in $L^2(0,T;H^{-1}(\T))$. Using the identities
$$
m^k=\frac{\nu\rho^k-a^k\rho^k}{\nu-1},
\quad
n^k=\frac{a^k\rho^k-\rho^k}{\nu-1},
$$
together with \eqref{eq:flux-decomp-m}, we deduce
\begin{equation}\label{eq:mandn}
\p_t m^k-\p_x J_m^k=\frac{\nu r_0^k-r_1^k}{\nu-1},
\quad
\p_t n^k-\nu\p_x J_n^k=\frac{r_1^k-r_0^k}{\nu-1}
\end{equation}
in $\mathcal D'((0,T)\times\T)$. The right-hand sides converge to $0$ in distributions. Since
$$
m^k\weakto m,
\quad
n^k\weakto n
\quad \text{weakly in } L^2((0,T)\times\T),
$$
and, by Lemma~\ref{lem:flux-identification},
$$
J_m^k\weakto m\,\p_x\rho,
\quad
J_n^k\weakto n\,\p_x\rho
\quad \text{weakly in } L^1((0,T)\times\T),
$$
we may pass to the limit in the previous distributional identities and obtain
$$
\p_t m-\p_x(m\,\p_x\rho)=0,
\quad
\p_t n-\nu\p_x(n\,\p_x\rho)=0
$$
in $\mathcal D'((0,T)\times\T)$. Also, using~\eqref{eq:mandn} we have that $m^k$ and $n^k$ are bounded uniformly in $H^{1}(0,T; H^{-1}(\T))$. Therefore by integrating by parts and using Assumption~\ref{A6} we obtain the initial traces $t=0$. Therefore $(m,n)$ is a weak solution in the sense of Definition~\ref{def:weak-solution}.
\end{proof}

\bibliographystyle{siam}
\bibliography{biblio}

\end{document}